\documentclass[a4paper,11pt]{article}
\usepackage{graphicx}
\usepackage{amsfonts}
\usepackage{amsmath}
\usepackage{amssymb}
\usepackage{indentfirst}
\usepackage{titlesec}

\def\barr{\begin{array}}
\def\earr{\end{array}}
\def\bali{\begin{aligned}}
\def\eali{\end{aligned}}
\def\bearr{\begin{eqnarray}}
\def\eearr{\end{eqnarray}}

\providecommand{\play}{\displaystyle}

\providecommand{\li}{\limits}

\providecommand{\pt}{\partial}

\providecommand{\ra}{\rightarrow}

\providecommand{\da}{\downarrow}

\providecommand{\Prob}{\mathbb{P}}

\providecommand{\E}{\mathbb{E}}

\providecommand{\al}{\alpha}

\providecommand{\bt}{\beta}

\providecommand{\gm}{\gamma}

\providecommand{\Gm}{\Gamma}

\providecommand{\dt}{\delta}

\providecommand{\Dt}{\Delta}

\providecommand{\ve}{\varepsilon}

\providecommand{\kp}{\kappa}

\providecommand{\lb}{\lambda}

\providecommand{\zt}{\zeta}

\providecommand{\N}{\mathbb N}

\providecommand{\R}{\mathbb R}

\providecommand{\cF}{\mathcal F}

\providecommand{\cH}{\mathcal H}

\providecommand{\cL}{\mathcal L}

\providecommand{\vphi}{\varphi}

\providecommand{\1}{\mathbf 1}

\providecommand{\hH}{\widehat{H}}

\providecommand{\Hbar}{\overline{H}}

\begin{document}

\title{On metastability in nearly-elastic systems}
\author{Wenqing Hu\thanks{Department of Mathematics,
University of Maryland at College Park, huwenqing@math.umd.edu}}
\date{}

\maketitle

\begin{abstract}
We consider a nearly-elastic model system with one degree of
freedom. In each collision with the "wall", the system can either
lose or gain a small amount of energy due to stochastic
perturbation. The weak limit of the corresponding slow motion, which
is a stochastic process on a graph, is calculated. A large deviation
type asymptotics and the metastability of the system is also
considered.
\end{abstract}



\section{Introduction}

Consider a model of a one-dimensional system with several potential
wells (Fig.1). A particle of unit mass moves freely in an interval
$[q_1,q_n]$ with elastic reflection at the ends of the interval if
the initial velocity is large enough. Let a finite number of points
$q_2, q_3,...,q_{n-1}\in (q_1, q_n)$ be given. Suppose at each $q_i$
there is a "wall" of certain height which gives the particle
instantaneous reflection once the particle hits it from either side.
The "height" coordinate $H$ is the energy of the particle. The
potential wells are numbered by $1,2,...,N$ (see Fig.1, where
$N=7$). Note that some of the wells are the combination of "smaller"
wells. For example, in Fig.1 well 5 consists of wells 1 and 2, well
6 consists of wells 5 and 3, and well 7 consists of wells 6 and 4.
The speed of the particle at energy level $H$ is $\sqrt{2H}$. In the
following, we always make the convention that the bigger wells, like
well 5 which consists of wells 1 and 2, are of energy level between
the top of that well and the one that separates the two smaller
wells. For example, in Fig.1 well 5 is supposed to be of energy
level between $H_6$ and $H_5$; well 6 is supposed to be of energy
level between $H_7$ and $H_6$, etc. Under this convention each well
with number $i$ has a minimum energy level $H_i$ (see Fig.1). We
assume that all $H_i$'s are bounded away from $0$. Within well $i$,
at energy level $H$, the particle moves between the walls of that
well and has speed $v=\sqrt{2H}$. At each collision with the wall,
the particle is instantaneously reflected and the speed of the
particle remains the same. The energy $H$ is preserved in the
system.

\begin{figure}
\centering
\includegraphics[height=7cm, width=15cm , bb=6 63 560 309]{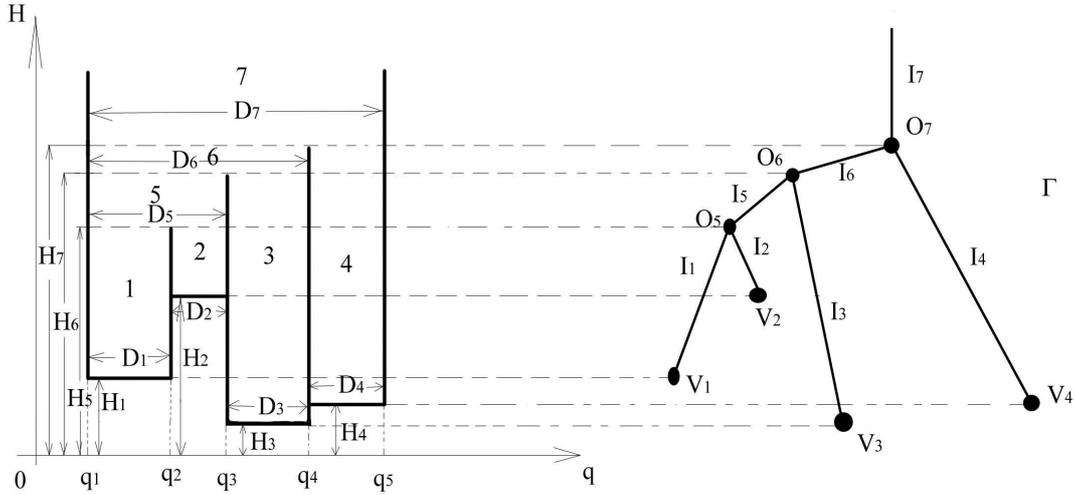}
\caption{The 1-dimensional mechanical model}
\end{figure}

Assume now that the collisions with the walls are not absolutely
elastic. If the particle is in well $i$ with energy $H$, then it
hits the left (right) wall of that well and was reflected, while at
the same time its energy becomes $H-\ve \xi^{(i)}_k$ ($H-\ve
\eta^{(i)}_k$) (if at the bottom of the well $i$ there is no smaller
wells the energy decays to $H_i\vee(H-\ve \xi^{(i)}_k)$ or $H_i\vee
(H-\ve \eta^{(i)}_k)$, respectively, and $a\vee b=\max (a,b)$). Here
$0<\ve<<1$ is a small parameter and $k$ denotes the number of
collisions with the left (right) wall (when the particle is at some
energy level which is the bottom of a "big" well, i.e., one which
contains two smaller wells we take $\xi_{k}^{(i)}$ and
$\eta_k^{(i)}$ to be those corresponding to the big well). The
sequences of random variables $\{\xi^{(i)}_k\}_{k \geq 1}$,
$\{\eta^{(i)}_k\}_{k \geq 1}$ are i.i.d. and mutually independent
with $\mathbb{E}(\xi^{(i)}_k+\eta^{(i)}_k)>0$. We assume that these
random variables are bounded $\Prob\{|\xi_k^{(i)}|\leq
M\}=\Prob\{|\eta_k^{(i)}|\leq M\}=1$ for some $M>0$ and they all
have continuous densities. In all the following, when we use random
variables such as $\xi$, $\eta$ without subscript, they are
understood as independent random variables and having the same
distribution as corresponding $\xi_k$ and $\eta_k$'s. Also, later in
this paper we will always denote $\zt_k=-(\xi_k+\eta_k)$ and
$\zt=-(\xi+\eta)$.

The position of the particle in our perturbed system can now be
described by a stochastic process
$\widetilde{X}_t^\ve=(\widetilde{H}_t^\ve, \widetilde{q}_t^\ve)$
where $\widetilde{H}_t^{\ve}$ is the energy of the particle at time
$t$ and $\widetilde{q}_t^\ve$ is the horizontal position of the
particle (see Fig.1). We denote the width of the $i$-th well by
$D_i$. In Fig.1 we have $D_5=D_1+D_2$, $D_6=D_5+D_3$ and
$D_7=D_6+D_4$.

The perturbed system $\widetilde{X}_t^\ve$ has, for $0<\ve<<1$, fast
and slow components. The fast component consists of the motion along
the non-perturbed trajectory. To describe the slow component,
consider the graph $\Gamma$ obtained after identification of points
of each well with a given energy level $H$. Denote by $\sqcap$ the
phase space of our system: $\sqcap$ is the union of all wells and it
is assumed that each interior well consists of two sides, left and
right. Denote by $Y: \sqcap \ra \Gamma$ the identification map of
the phase space $\sqcap$ to $\Gamma$. The slow component of the
motion is $\widetilde{Y}_t^\ve=Y(\widetilde{H}_t^\ve,
\widetilde{q}_t^\ve)$ (compare with [5, Ch.8], [4]). We rescale time
$t \mapsto t/\ve$. Define $X_t^\ve=\widetilde{X}_{t/\ve}^\ve$,
$H_t^\ve=\widetilde{H}_{t/\ve}^\ve$,
$q_t^\ve=\widetilde{q}_{t/\ve}^\ve$,
$Y_t^\ve=\widetilde{Y}_{t/\ve}^\ve$.

We make a convention here: in the following processes with a tilde
on it are original processes with natural time parameter $t$;
processes without such a tilde on it are time-rescaled process with
time $t/\ve$; processes with a hat on it are piecewise linear
modifications of the one under the hat. For example,
$H^\ve_t=\widetilde{H}^\ve_{t/\ve}$ and
$\widehat{\widetilde{H}^\ve_t}$ is a piecewise linear modification
of $\widetilde{H}_t^\ve$, $\widehat{H}^\ve_t$ is a piecewise linear
modification of $H_t^\ve$, etc. Here piecewise linear modifications
are obtained by joining each consecutive corners of the step
functions $\widetilde{H}_t^\ve$ and $H_t^\ve$.

Number the edges of the graph: $\Gamma=\{I_1, I_2,..., I_N\}$ where
$N$ is the number of the wells (in Fig.1 $N=7$). The $i$-th well
corresponds to edge $I_i$. Exterior vertex $V_k$ corresponds to the
bottom of the $k$-th well. Interior vertex $O_l$ corresponds to the
lowest energy level (as was in the convention made before) of the
$l$-th well ("big" well). Then $Y(H,q)=(H,K(H,q))$ where $K(H,q)$ is
the number of the edge containing $Y(H,q)$ and $H$ is the energy. We
see that after time rescaling, the slow component is the process
$Y_t^\ve=(H_t^\ve, K(H_t^\ve, q_t^\ve))$.

We will show that the process $Y_t^\ve$ converges, as $\ve \da 0$,
to a stochastic process $Y_t$ on $\Gamma$. The process $Y_t$ is a
deterministic motion within each edge of $\Gamma$ and has
stochasticity only at the interior vertices $O_l$ of $\Gamma$.

Since we allow random variables $\xi_k^{(i)}$, $\eta_k^{(i)}$ to be
less than 0, it can happen that the particle enters certain well and
sooner or later it jumps out of that well. Since we assumed that
$\mathbb{E}(\xi_k^{(i)}+\eta_k^{(i)})>0$, this is a large deviation
type event. We will calculate the "quasi-potential" describing how
difficult it is to switch from one well to another. For the system
with many wells, metastability and asymptotic behavior of the system
will be considered in Section 4.

\section{The limiting process}

In this section we first consider the two well case. Let us assume
that our system has two wells 1 and 2 and their combination is well
3. Interior vertex is $O_3$ and exterior vertices are $V_1$ and
$V_2$. The edges are $I_1$, $I_2$ and $I_3$. We assume that
$\mathbb{P}\{|\xi_k^{(i)}|\leq M\}=\mathbb{P}\{|\eta_k^{(i)}|\leq
M\}=1$ for some constant $M>0$. Similar to [4], by using the
standard averaging principle, we get

\

\textbf{Lemma 2.1.} \textit{Let $H^\ve_0=H_0>H(O_3)$. Within each
edge of the graph $\Gamma$, as $\ve \da 0$, the process
$H^\ve_t=\widetilde{H}^\ve_{t/\ve}$, converges uniformly in
probability on $0<t<T<\infty$, to a deterministic motion $H(t)$
which is defined by the equations}

$$H(t)=\left(\sqrt{H_0}-t\play{\frac{\mathbb{E}\xi^{(3)}+\mathbb{E}\eta^{(3)}}{2\sqrt{2}D_3}}\right)^2, 0<t\leq t_0 \
on \ I_3; \eqno(2.1)$$ \textit{and}
$$H(t)=\left(\sqrt{H(O_3)}-(t-t_0)\play{\frac{\mathbb{E}\xi^{(1)}+\mathbb{E}\eta^{(1)}}{2\sqrt{2}D_1}}\right)^2, t>t_0 \
on \ I_1; \eqno(2.2)$$
$$H(t)=\left(\sqrt{H(O_3)}-(t-t_0)\play{\frac{\mathbb{E}\xi^{(2)}+\mathbb{E}\eta^{(2)}}{2\sqrt{2}D_2}}\right)^2, t>t_0 \
on \ I_2 \eqno(2.3)$$ \textit{respectively. Here $H(O_3)$ is the
energy corresponding to the interior vertex $O_3$ and
$t_0=\play{\frac{2\sqrt{2}D_3(\sqrt{H_0}-\sqrt{H(O_3)})}{\mathbb{E}\xi^{(3)}+\mathbb{E}\eta^{(3)}}}$
is the time for $H(t)$ to come to the interior vertex $O_3$.}

\

Similarly as was done in [4], we consider a piecewise linear
modification $\widehat{\widetilde{H}}^\ve_t$ of
$\widetilde{H}^\ve_t$. Under the convention made in the introduction
we put $\widehat{H}_t^\ve=\widehat{\widetilde{H}}_{t/\ve}^\ve$ and
$\widehat{X}_t^\ve=(\widehat{H}_t^\ve, q_t^\ve)$. Let
$\widehat{Y}_t^\ve=(\widehat{H}^\ve_t, K(\widehat{X}_t^\ve))$. It is
clear that for fixed $\ve >0$,
$$\mathbb{P}\{|\widehat{\widetilde{H}}^\ve_t-\widetilde{H}^\ve_t|<C\ve\}=1 \eqno(2.4)$$ for some
positive constant $C>0$ and $0<t<T<\infty$. We have, as in [4],

\

\textbf{Lemma 2.2.} \textit{For each $T>0$, the family
$\{\widehat{Y}^\ve_t\}_{t>0}$ is tight in $C_{0T}(\Gamma)$. }

\

We now turn to the problem of determining the asymptotic branching
probability for the process $\widehat{Y}_t^\ve$ as $\ve \da 0$, at
the interior vertex $O_3$. Let us first present an auxiliary lemma
about certain properties of random walk (compare with [4]).

Let $\{\xi_k\}_{k\geq 1}$, $\{\eta_k\}_{k \geq 1}$ be i.i.d,
mutually independent sequences of random variables. Assume that the
random variables have continuous densities and
$\mathbb{P}\{-\infty<-\al<\xi_k<\al<\infty\}=1$,
$\mathbb{P}\{0<\eta_k<\al<\infty\}=1$ for some positive constant
$\al>0$. Notice that we allow $\xi_k$ to be negative but we assume
that $\mathbb{E}(\xi_k+\eta_k)>0$. Put, for $m\geq 0$, that

$$S_0=0 \ , \ S_{2m}=\sum_{k=1}^{m}(\xi_k+\eta_k) \ , \ S_{2m+1}=S_{2m}+\xi_{m+1} \ .$$

Define $\tau_{n}^\lambda=\min\{m: S_m>n\lambda\}$ for $\lambda
>0$.

Since $\mathbb{E}(\xi_k+\eta_k)>0$, the law of large numbers implies
that $\mathbb{P}\{\tau_n^\lambda<\infty\}=1$ for any $\lambda>0$, $n
\in \mathbb{Z}$.

Let $\zeta_k=\xi_k+\eta_k$. Let $T_n=S_{2n}=\sum\li_{k=1}^n
\zeta_k$. Sample trajectories of $S_n$ and $T_n$ are shown in Fig.2.
Let $\mathbf{N}=\min\li_{k \geq
 1}\{T_n>0\}$. Put

$$E_n(I)=\mathbb{P}\{\mathbf{N}=n, T_n\in I\}$$
for $I \subset (0, +\infty)$. In other words, $E_n(I)$ is the
probability of the event
$$\{T_1\leq 0, T_2 \leq 0, ... , T_{n-1} \leq 0, T_n>0, T_n\in I\}.$$

Consider random variables $a=T_{\mathbf{N}-1}$,
$b=S_{2\mathbf{N}-1}$, $c=T_{\mathbf{N}}$. We are now ready to state

\

\begin{figure}
\centering
\includegraphics[height=7cm, width=15cm , bb=6 63 560 309]{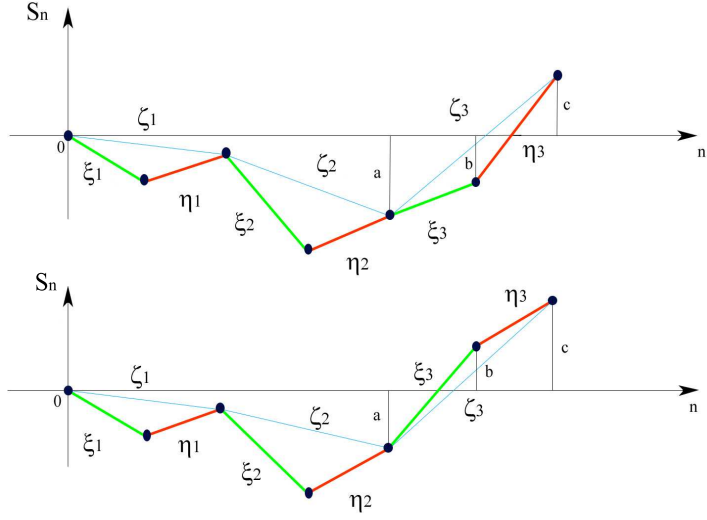}
\caption{Strong ascending ladder points}
\end{figure}

\

\textbf{Lemma 2.3.} \textit{Under mentioned above conditions,}

$$\lim\li_{n \ra \infty}\mathbb{P}\{\tau_n^\lambda \ is \ odd \}=\play{\frac{\mathbb{E}b\mathbf{1}_{b>0}}{\mathbb{E}c}} \ ,$$
$$\lim\li_{n \ra \infty}\mathbb{P}\{\tau_n^\lambda \ is \ even \}=\play{1-\frac{\mathbb{E}b\mathbf{1}_{b>0}}{\mathbb{E}c}} \ .$$

\

\textbf{Proof.} We say a \textit{strong ascending ladder point} (see
[3, Ch.12]) for $\{T_n\}_{n \geq 1}$ (respectively, $\{S_n\}_{n \geq
1}$) occurs at step $k$ if
$$T_k>\max\{T_r: 0 \leq r \leq k-1\}$$
(respectively, $S_k>\max\{S_r: 0 \leq r \leq k-1\}$).

If the successive strong ascending ladder points for $T_n$ are
$W_1$, $W_1+W_2$, ... , we write
$T_{W_1+W_2+...+W_k}-T_{W_1+W_2+...+W_{k-1}}$ (suppose $W_0=0$) as
$Z_k$, the $k$-th \textit{strong ascending ladder step} for
$\{T_n\}_{n \geq 1}$.

The random variables $W_k$, $k \geq 1$ are i.i.d with common
distribution the same as that of $\mathbf{N}$. The random variables
$Z_k$, $k\geq 1$ are i.i.d with common distribution the same as that
of $c=T_{\mathbf{N}}$.

Since we assumed that $\mathbb{P}\{0<\eta_k<\al<\infty\}=1$, the
occurrence of a strong ascending ladder point for $\{T_n\}_{n \geq
1}$ at step $k$ implies that a strong ascending ladder point for
$\{S_n\}_{n \geq 1}$ happens either at step $2k$ or at step $2k-1$
(see Fig.2). Define
$R_k=S_{2(W_1+W_2+...+W_k)-1}-S_{2(W_1+W_2+...+W_{k-1})}$. Since
each piece of the random walk between steps $W_1+...+W_{k-1}$ and
$W_1+...+W_k$ are i.i.d, the random variables $R_k$ are i.i.d with
common distribution the same as that of $b=S_{2\mathbf{N}-1}$. By
using the same local limit argument as that in [4, Lemma 3.3], one
can see that

$$\lim\li_{n \ra \infty}\mathbb{P}\{\tau_n^\lambda \ \text {is \
odd} \}=\play{\lim\li_{n \ra \infty}\frac{\sum\li_{k=1}^{n}
R_k\mathbf{1}_{R_k>0}
}{\sum\li_{k=1}^{n}Z_k}}=\play{\frac{\mathbb{E}b\mathbf{1}_{b>0}}{\mathbb{E}c}}
\ ,$$ and the result follows. $\square$

\

By using the Wiener-Hopf theory (see [3, Ch.12, Ch.18]), one can
sometimes determine the distribution of $c$ by
$E(dx)=\sum\li_{n=0}^\infty E_n(dx)$ and thus
$\mathbb{E}c=\play{\int_0^\infty xE(dx)}$. The distribution of $b$
can be determined by the convolution relation
$$E(dx)=\int_{-\infty}^x b(dy)F_{\eta}(dx-y).$$

Here $b(dy)=\mathbb{P}\{b\in dy\}$ and $F_\eta$ is the common
distribution function of $\eta_k$. We refer the reader to [3, Ch.12,
Ch.18].

\

Now let us turn back to our system. Assume that our system always
loses energy on the right walls, i.e.
$\mathbb{P}\{0<\eta_k^{(i)}<M<\infty\}=1$. On the left walls the
system can either gain or lose energy - we only assume that
$\mathbb{P}\{-\infty<-M<\xi_k^{(i)}<M<\infty\}=1$. Let $a^{(i)},
b^{(i)}, c^{(i)}$ be defined in the same way as $a,b,c$ in Lemma 2.3
for random walks $S_n^{(i)}$ constructed from $\xi^{(i)}_k,
\eta^{(i)}_k$:

$$S_0^{(i)}=0 \ , \ S_{2m}^{(i)}=\sum_{k=1}^{m}(\xi_k^{(i)}+\eta_k^{(i)}) \ , \ S_{2m+1}^{(i)}=S_{2m}^{(i)}+\xi_{m+1}^{(i)} \ .$$

Let $\lambda=H_0-H(O_3)$, $n=\left[\play{\frac{1}{\ve}}\right]$, and
apply Lemma 2.3 directly, we get

\

\textbf{Lemma 2.4.}

$\lim\li_{\ve \da 0}\mathbb{P}\{X_t^\ve \ finally \ falls \ into \
well \ 1\} =
\play{\frac{\mathbb{E}b^{(3)}\mathbf{1}_{b^{(3)}>0}}{\mathbb{E}c^{(3)}}}
\equiv p^{(3)}_1,$

$\lim\li_{\ve \da 0}\mathbb{P}\{X_t^\ve \ finally \ falls \ into \
well \ 2\} =
\play{1-\frac{\mathbb{E}b^{(3)}\mathbf{1}_{b^{(3)}>0}}{\mathbb{E}c^{(3)}}}\equiv
p^{(3)}_2.$

\

Define a process $Y_t$ on $\Gamma$: $Y_t=(H(t), K(t))$; on each edge
$H(t)$ satisfies equations (2.1), (2.2), (2.3), respectively. The
process $Y_t$, when arriving at the interior vertex $O_3$,
immediately leaves that vertex and goes into edge $I_1$ or $I_2$
with probabilities $p^{(3)}_1$ and $p^{(3)}_2$, respectively. We
have,

\

\textbf{Theorem 2.1.} \textit{Under the same assumption mentioned
before Lemma 2.4, as $\ve \da 0$, process $\widehat{Y}^\ve_t$
converges weakly, for $0<T<\infty$, in $C_{0T}(\Gamma)$ with uniform
topology, to $Y_t$. }

\

The above result can be easily generalized to the case when the
system has more than two wells. The averaging principle is the same
as before: within each edge $I_i$, as $\ve \da 0$, $H^\ve_t$
converges to a deterministic motion $H(t)$ which satisfies the
differential equation

$$\frac{dH}{dt}=-\frac{\mathbb{E}\xi^{(i)}+\mathbb{E}\eta^{(i)}}{T_i(H)} \ ,$$
where $i$ is the number of the well,
$T_i(H)=\play{\frac{2D_i}{\sqrt{2H}}}$ is the period of the elastic
motion within well $i$.

We assume that the system always loses energy on the right walls,
i.e. $\mathbb{P}\{0<\eta_k^{(i)}<M<\infty\}=1$; and on the left
walls the system can either gain or lose energy - we only assume
that $\mathbb{P}\{-\infty<-M<\xi_k^{(i)}<M<\infty\}=1$. The
branching probabilities for the limiting motion $Y_t$ at the bottom
of well $i$ can be given by
$p_1^{(i)}=\play{\frac{\mathbb{E}b^{(i)}\mathbf{1}_{b^{(i)}>0}}{\mathbb{E}c^{(i)}}}$
(for entering the left well) and
$p_2^{(i)}=\play{1-\frac{\mathbb{E}b^{(i)}\mathbf{1}_{b^{(i)}>0}}{\mathbb{E}c^{(i)}}}$
(for entering the right well). The branching at each interior vertex
is independent of the others.

\

Finally we briefly consider the case when we throw away the
artificial restriction that $\Prob\{0<\eta_k<\al<\infty\}=1$.
Suppose $\xi_k$ and $\eta_k$ are two i.i.d. series and mutually
independent. Let $\Prob\{-\al< \xi_k < \al\}=\Prob\{-\al< \eta_k<
\al\}=1$ for some $\al>0$ and $k=1,2,3, ... $ . Suppose all
$\xi_k$'s and $\eta_k$'s have continuous densities. We also assume
that $\E(\xi_k+\eta_k)>0$. Let us add one more assumption that
$\Prob\{\xi_k>0\}>0$, $\Prob\{\eta_k>0\}>0$.

Let us consider the strong ascending ladder points for the random
walk

$$S_0=0 \ , \ S_{2m}=\sum\li_{k=1}^m (\xi_k+\eta_k) \ , \ S_{2m+1}=S_{2m}+\xi_{m+1} \ .$$

We define these strong ascending ladder points to be $J_1$,
$J_1+J_2$, ... . Let $M_k=1$ if $J_1+...+J_k$ is odd and $M_k=0$ if
$J_1+...+J_k$ is even. Let $M_0=0$.

Let us consider another random walk

$$S_0'=0 \ , \ S_{2m}'=\sum\li_{k=1}^m (\eta_k+\xi_k) \ , \ S_{2m+1}'=S_{2m}'+\eta_{m+1} \ $$
and the corresponding strong ascending ladder points $J_1'$,
$J_1'+J_2'$, ... . We consider the first strong ascending ladder
steps $\gamma^{(0)}=S_{J_1}$ for $\{S_n\}_{n\geq 0}$ and
$\gamma^{(1)}=S'_{J_1'}$ for $\{S_n'\}_{n\geq 0}$.

By strong Markov property of the random walk $S_n$ and our
assumptions on $\xi_k$ and $\eta_k$ it is easy to see that $M_k$ is
an ergodic Markov chain with two states $\{0,1\}$ and an invariant
measure $\mu(\{0\})=\mu_0$ and $\mu(\{1\})=\mu_1$ for some
$0<\mu_i<1$ and $\sum\mu_i=1$, $i=0,1$. The coupling chain
$(M_{k-1}, M_k)$ is also an ergodic Markov chain with four states
$\{(0,0), (0,1), (1,0), (1,1)\}$ and an invariant measure
$\mu\{(0,0)\}=\mu_{00}$, $\mu\{(0,1)\}=\mu_{01}$,
$\mu\{(1,0)\}=\mu_{10}$, $\mu\{(1,1)\}=\mu_{11}$. Here
$0<\mu_{ij}<1$ and $\sum\mu_{ij}=1$ for $i,j=0,1$. It is clear that
$\mu_{11}+\mu_{10}=\mu_1$, $\mu_{01}+\mu_{00}=\mu_0$.

Let $\gm_k^{(0)}$ be a sequence of i.i.d random variables which has
common distribution same as $\gm^{(0)}$. Let $\gm_k^{(1)}$ be a
sequence of i.i.d random variables which has common distribution
same as $\gm^{(1)}$. The random variables $\gm^{(0)}$ and
$\gm^{(1)}$ are bounded and have continuous densities. We choose
these random variables such that they are mutually independent and
also independent of the $M_k$'s.

Define $\tau_{n}^\lambda=\min\{m: S_m>n\lambda\}$ for $\lambda
>0$.

We claim the

\

\textbf{Lemma 2.5.} \textit{Under mentioned above conditions,}

$$\lim\li_{n \ra \infty}\mathbb{P}\{\tau_n^\lambda \ is \ odd \}
=\dfrac{\mu_{11}\E\gm^{(1)}+\mu_{01}\E\gm^{(0)}}{\mu_1\E\gm^{(1)}+\mu_0\E\gm^{(0)}}
\ ,$$
$$\lim\li_{n \ra \infty}\mathbb{P}\{\tau_n^\lambda \ is \ even \}
=\dfrac{\mu_{10}\E\gm^{(1)}+\mu_{00}\E\gm^{(0)}}{\mu_1\E\gm^{(1)}+\mu_0\E\gm^{(0)}}
\ .$$

\

\textbf{Proof.} We use the same local limit theorem argument as in
[4, Lemma 3.3]. We first apply the local limit theorem to sequence
$T_n$ (as defined in the proof of Lemma 2.3). Then we use the fact
that

$$\begin{array}{l}
\lim\li_{n\ra \infty}\Prob\{\tau_n^\lb\text{ is odd }\}
\\
=\lim\li_{n \ra \infty}\dfrac{\sum\li_{k=1}^n M_k
\gm_k^{(M_{k-1})}}{\sum\li_{k=1}^n \gm_k^{(M_{k-1})}}
\\
=\lim\li_{n \ra
\infty}\dfrac{\dfrac{\nu_{11}(n)}{n}\dfrac{1}{\nu_{11}(n)}\sum\li_{k=1}^{\nu_{11}(n)}
\gm_k^{(1)}+\dfrac{\nu_{01}(n)}{n}\dfrac{1}{\nu_{01}(n)}\sum\li_{k=1}^{\nu_{01}(n)}
\gm_k^{(0)}}{\dfrac{\nu_{1}(n)}{n}\dfrac{1}{\nu_{1}(n)}
\sum\li_{k=1}^{\nu_{1}(n)}
\gm_k^{(1)}+\dfrac{\nu_{0}(n)}{n}\dfrac{1}{\nu_{0}(n)}
\sum\li_{k=1}^{\nu_{0}(n)} \gm_k^{(0)}}
\\
=\dfrac{\mu_{11}\E\gm^{(1)}+\mu_{01}\E\gm^{(0)}}{\mu_1\E\gm^{(1)}+\mu_0\E\gm^{(0)}}
\
\end{array}$$
and the Lemma follows. Here for $i,j=0,1$ we set
$$\nu_{ij}(n)=\text{ number of } k  \text{'s such that } (M_{k-1},
M_{k})=(i,j), 1\leq k \leq n$$ and for $i=0,1$ we set
$$\nu_{i}(n)=\text{ number of } k \text{'s such that } M_{k-1}=i ,
1\leq k \leq n \ .$$ $\square$

\

But in this case the process $T_n$ loses its ability to "detect" a
strong ascending ladder point for the process $S_n$. Actually it
might happen that $S_1,...,S_{2n}$ have a strong ascending ladder
point at $S_{2n-1}$, yet $T_1,...,T_n$ have no strong ascending
ladder point. Therefore one might not get explicit formulas as in
Lemma 2.3. This problem of explicitly calculating the asymptotic
branching probability still remains open.

\

Now we turn back to our original system. By the same arguments that
we use to prove Theorem 2.1 we assert that under the assumptions
made in Section 1 and an additional assumption
$\Prob\{\xi_k^{(i)}>0\}>0$, $\Prob\{\eta_k^{(i)}>0\}>0$, we have

\

\textbf{Theorem 2.2.} \textit{As $\ve\da 0$ the process
$\widehat{Y}_t^\ve$ converges weakly for $0<T<\infty$ in
$C_{0T}(\Gm)$ with uniform topology to a process $Y_t$ on $\Gm$
which is a Markov process on $\Gm$. It is deterministic inside the
edges and only has stochasticity (i.e. certain branching
probabilities) at the interior vertices.}

\section{Large deviations}

We now calculate large deviation type asymptotics. We consider the
simplest case when there is only one well. The general case follows
from our result for one well case and will be discussed in the next
section.

Suppose our well has width $D$. The perturbation for the collision
at the walls is given by i.i.d and mutually independent sequences
$\{\xi_k\}_{k\geq 1}$ and $\{\eta_k\}_{k \geq 1}$. We assume that
$\mathbb{P}\{-M\leq\xi_k\leq M\}=\mathbb{P}\{-M\leq\eta_k\leq M\}=1$
for some $M>0$. We assume that $\mathbb{E}(\xi_k+\eta_k)>0$. Both
$\xi$ and $\eta$ have continuous density. This implies that the
process $H_t^\ve$ is bounded for time $0\leq t \leq T <\infty$. Let
us assume that for the time $0\leq t \leq T$ we have $0<H_0\leq
H^\ve_t\leq \Hbar<\infty$. Let $q_0^\ve=q_0$. Let
$$Q_t^\ve = q_0 + \text{the  total  horizontal  distance  that}  \,
q_t^\ve \,
 \text{traveled up to time} \, t \ .$$

The system $(H_t^\ve, Q_t^\ve)$, satisfies the equations:

$$
\left\{
     \begin{array}{l}
             \play{\dot{H}^\ve_t=-f(Q_t^\ve)} \ ,\\
             \play{\dot{Q}_t^\ve=\frac{1}{\ve}\sqrt{2H_t^\ve}} \ .
     \end{array}
\right.
 \eqno(3.1)
$$

Here random function $f(Q)=\sum\li_{k=1}^\infty\left(
\xi_k\delta(Q-(2k-1)D)+\eta_k\delta(Q-2kD)\right)$ where
$\delta(\cdot)$ is the Dirac $\delta$-function.

Consider a piecewise linear modification
$\widehat{\widetilde{H}}_t^\ve$ of the step function
$\widetilde{H}_t^\ve$, as defined at the beginning of Section 2. We
see that by (2.4) $\widehat{\widetilde{H}}_t^\ve$ is a good
approximation of $\widetilde{H}_t^\ve$.

System (3.1) has fast component $Q$ and slow component $H$ and they
depend on each other. Let $H_0\leq h\leq \Hbar$. Let
$Q^h(t)=q_0+t\sqrt{2h}$. Let $\beta \in \mathbb{R}$. Define

$$\begin{array}{l}
\mathcal{H}(h,\beta)
\\
\play{=\lim\li_{T\ra \infty}\frac{1}{T}\ln \mathbb{E}
\exp\left(-\beta\int_0^T f(Q^h(t))dt\right)}
\\
\play{=\dfrac{\sqrt{2h}}{2D}\ln\mathbb{E}\exp\left(-\beta(\xi+\eta)\right)}
\ .
\end{array}\eqno(3.2)$$

Let $\mathcal{L}$ be the Legende transform of $\mathcal{H}$:
$$\mathcal{L}(h,\alpha)=\sup\li_{\beta}(\alpha\beta-\mathcal{H}(h,\beta)).$$

Let $\vphi\in C_{[0,T]}([H_0,\Hbar])$. Let

$$S_{0T}(\varphi)=\left\{
\begin{array}{l}
\play{\int_0^T \mathcal{L}(\varphi_s, \dot{\varphi}_s)ds}  \ , \
\text{ for } \varphi \text{
absolutely continuous} , \\
+\infty \ ,  \text{ otherwise } \ .
\end{array}\right. \eqno(3.3)$$

We have:

\

\textbf{Theorem 3.1.} \textit{The family $\widehat{H}^\ve_t \ , \
0<t<T $ satisfies the large deviation principle as $\ve \da 0$ in
the space $C_{[0,T]}([H_0,\overline{H}])$ with normalizing factor
$\ve^{-1}$ and an action functional $S_{0T}(\vphi)$.}

\

To be precise, Theorem 3.1 means the following (see [5, Ch.3]):

(0) The set $\Phi(s)=\{\varphi\in C_{[0,T]}([H_0, \Hbar]):
S_{0T}(\varphi)\leq s\}$ is compact for every $s\geq 0$.

(I) For any $\nu>0$, any $\dt>0$, any $\vphi\in
C_{[0,T]}([H_0,\Hbar])$, there exist $\ve_0>0$ such that for any
$0<\ve\leq\ve_0$ we have

$$\Prob\{\rho_{0T}(\widehat{H}^\ve, \vphi)<\dt\}\geq
\exp(-\ve^{-1}(S_{0T}(\vphi)+\nu)) \, .$$

(II) For any $\dt>0$, any $\nu>0$ and any $s>0$ there exist an
$\ve_0>0$ such that for any $0<\ve\leq\ve_0$ we have

$$\Prob\{\rho_{0T}(\hH^\ve, \Phi(s))\geq \dt\}\leq\exp(-\ve^{-1}(s-\nu)) \, .$$

\

Here for $\vphi, \psi\in C_{0T}([H_0,\Hbar])$ we denote
$\rho_{0T}(\vphi, \psi)=\max\li_{0\leq t \leq T}|\vphi(t)-\psi(t)|$
and $\rho_{0T}(\vphi, \Phi(s))=\max\li_{\psi\in
\Phi(s)}\max\li_{0\leq t \leq T}|\vphi(t)-\psi(t)|$.

\

Let us consider an example.

\

\textbf{Example.} Let function $\cH_0(\beta)=\ln
\mathbb{E}\exp(-\beta(\xi+\eta))$. We have
$\mathcal{H}(h,\beta)=\play{\frac{\sqrt{2h}}{2D}}\cH_0(\beta)$.
Using the convexity of exponential function, we get
$\cH_0(\beta)\geq \ln \exp (-\beta \mathbb{E}(\xi+\eta))=-\beta
\mathbb{E} (\xi+\eta)$, i.e. $\cH_0(\beta)+\beta
\mathbb{E}(\xi+\eta) \geq 0$. The minimum is achieved at $\beta=0$.

Now we let $\varphi_t=H(t)$. Here $H(t)$ is the limiting motion of
$H^\ve_t$ as $\ve \da 0$. Standard averaging principle gives us

$$\frac{dH(t)}{dt}=-\frac{\sqrt{2H(t)}}{2D}(\mathbb{E}\xi+\mathbb{E}\eta) \ , \ H(0)=H_0^\ve \ .$$

Now $\mathcal{H}(H(t),
\beta)=\play{\frac{\sqrt{2H(t)}}{2D}}\cH_0(\beta)$, and

$$\begin{array}{l}
\mathcal{L}(H(t),\dot{H}(t))\\
=\sup\li_{\beta}(\dot{H}(t)\beta-\mathcal{H}(H(t),\beta))\\
=\sup\li_{\beta}(\dot{H}(t)\beta-\play{\frac{\sqrt{2H(t)}}{2D}}\cH_0(\beta))\\
=-\play{\frac{\sqrt{2H(t)}}{2D}}
\inf\li_\beta((\mathbb{E}\xi+\mathbb{E}\eta)\beta+\cH_0(\beta))=0  \
.
\end{array}$$

This means that $S_{0T}(H(t))=0$ , which is not surprising since
$H(t)$ is the averaged motion of the system.

On the other hand, for any absolutely continuous trajectory
$\vphi_t$ such that $\dot{\vphi}_t>0$ for $0\leq t \leq T$ we have
$\cL(\vphi_t,\dot{\vphi}_t)=\sup\li_{\bt}(\dot{\vphi}_t \bt -
\cH(\vphi_t, \bt))>0$ since $\cH(\vphi_t,0)=0$ and $\dfrac{\pt}{\pt
\bt}\cH(\vphi_t,\bt)|_{\bt=0}=-\dfrac{\sqrt{2\vphi_t}}{2D}\E(\xi+\eta)<0$
(see Lemma 3.2.1). This gives $S_{0T}(\vphi)>0$ which means that
there is a "difficulty" for the system to gain some energy.
$\square$

\

The \textbf{Proof} of Theorem 3.1 is based on a combination of
Cram\'{e}r's large deviation principle for i.i.d. sums and the
technique to calculate large deviations from an averaged system with
full dependence, which was developed in [6], [7].

Let us first formulate an analogue of the classical Cram\'{e}r's
large deviation principle for i.i.d. sums (compare with, for
example, [5, Ch. 5, Examples 1.3 and 1.4]). Our proof follows [2,
Section 2.2].

\

\textbf{Lemma 3.1.} \textit{Let $\zt_1$, ..., $\zt_n$,... be a
sequence of bounded i.i.d random variables, which have continuous
densities and let
$$\cH_0(\bt)=\ln \E \exp(\bt \zeta_i) \ .$$ Let
$\cL_0(\al)=\sup\li_{\bt\in \R}(\al\bt-\cH_0(\bt))$. Let $\ve>0$.
Suppose integer $n(\ve)\ra \infty$ as $\ve \da 0$. Suppose for a
bounded $-\infty<\underline{x}<x(\ve)<\overline{x}<\infty$ we have
$|\text{argmax}_{\bt}(x(\ve)\bt-\cH_0(\bt))|\leq b<\infty$ is
uniformly bounded. Then for any $\nu>0$, there exist
$\overline{\dt}>0$ such that, for any
$0<\underline{\dt}<\overline{\dt}$ and for any
$0<\underline{\dt}<\dt(\ve)<\overline{\dt}<\infty$, there exist
$\ve_0>0$ such that for any $0<\ve<\ve_0$ we have}

$$\exp(-n(\ve) (\cL_0(x(\ve))+\nu))\leq\Prob\left\{\left|\dfrac{\zt_1+...+\zt_{n(\ve)}}{n(\ve)}-x(\ve)\right|
<\dt(\ve)\right\}\leq \exp(-n(\ve) (\cL_0(x(\ve))-\nu)) \ .$$

\

\textbf{Proof.} Let $A(n)=\dfrac{1}{n}(\zt_1+...+\zt_n)$. Let $A=\E
A(n)$. We are estimating
$$\Prob\{x(\ve)-\dt(\ve)<A(n(\ve))<x(\ve)+\dt(\ve)\} \ .$$

About the upper bound. Consider first the case $x(\ve)-\dt(\ve)>A$.
We have, using Chebyshev inequality, for $\bt \geq 0$, that
$$\begin{array}{l}
\Prob\{x(\ve)-\dt(\ve)<A(n(\ve))\}
\\
\leq \exp(-\bt n(\ve)(x(\ve)-\dt(\ve)))\E(\bt n(\ve)A(n(\ve)))
\\
=\exp(-\bt n(\ve)(x(\ve)-\dt(\ve)))\E\prod\li_{k=1}^{n(\ve)}\exp(\bt
\zt_k)
\\
=\exp(-n(\ve)((x(\ve)-\dt(\ve))\bt-\cH_0(\bt))) \ .
\end{array}$$

Since for $x(\ve)-\dt(\ve)>A$ and $\bt \geq 0$ we have
$\cL_0(x(\ve)-\dt(\ve))=\sup\li_{\bt \geq 0}((x(\ve)-\dt(\ve))\bt -
\cH_0(\bt))$, we optimize the above inequality and we get
$$\Prob\{x(\ve)-\dt(\ve)<A(n(\ve))\}\leq \exp(-n(\ve)\cL_0(x(\ve)-\dt(\ve))) \ .$$

Since our choice of $x(\ve)$ makes
$|\text{argmax}_\bt(x(\ve)\bt-\cH_0(\bt))|$ uniformly bounded, the
uniform continuity of $\cL_0$ gives the upper bound in this case.
That is, we can choose $\overline{\dt}>0$ small enough such that for
$0<\underline{\dt}<\dt(\ve)<\overline{\dt}$ we have
$\cL_0(x(\ve)-\dt(\ve))\geq \cL_0(x(\ve))-\nu$.

In the case when $x(\ve)+\dt(\ve)<A$, we estimate, for $\bt \geq 0$,
that

$$\begin{array}{l}
\Prob\{-(x(\ve)+\dt(\ve))<-A(n(\ve))\}
\\
\leq \exp(\bt n(\ve)(x(\ve)+\dt(\ve)))\E(-\bt n(\ve)A(n(\ve)))
\\
=\exp(\bt n(\ve)(x(\ve)+\dt(\ve)))\E\prod\li_{k=1}^{n(\ve)}\exp(-\bt
\zt_k)
\\
=\exp(-n(\ve)((x(\ve)+\dt(\ve))(-\bt)-\cH_0(-\bt))) \ .
\end{array}$$

Now we use the fact that for $x(\ve)+\dt(\ve)<A$ we have
$\cL_0(x(\ve)+\dt(\ve))=\sup\li_{\bt \geq 0}((x(\ve)+\dt(\ve))(-\bt)
- \cH_0(-\bt))$ and we apply a similar argument.

Now in the case of $x(\ve)-\dt(\ve)\leq A\leq x(\ve)+\dt(\ve)$, we
choose $\overline{\dt}>0$ small enough such that
$|\cL_0(x(\ve))-\cL_0(A)|<\nu/2$ and we notice that $\cL_0(A)=0$.
This gives the trivial upper bound as $\ve \da 0$.

Now we prove the lower bound. Consider the unique solution of the
equation $$\cH_0'(\eta)=x(\ve) \ .$$ By our assumptions on the
uniform boundedness of $|\text{argmax}_\bt(x(\ve)\bt-\cH_0(\bt))|$
and about the boundedness an having density of $\zt$'s it is easy to
check that the solution of this equation exists and is unique. Now
define a new measure $\widehat{\Prob}^\ve$ in terms of $\Prob$ as
$$\dfrac{d\widehat{\Prob}^\ve}{d\Prob}(x)=\exp(\eta x-\cH_0(\eta)) \ .$$
This $\widehat{\Prob}^\ve$ is a probability measure since
$\play{\int_{\R}d\widehat{\Prob}^\ve=\dfrac{1}{\E
\exp(\eta\zt)}\int_{\R}\exp(\eta x)d\Prob}=1$. Also under
$\widehat{\Prob}^\ve$ the expected value of $\zt$ is
$\widehat{\E}^\ve\zt=\dfrac{1}{\E
\exp(\eta\zt)}\play{\int_{\R}x\exp(\eta
x)d\Prob}=\cH_0'(\eta)=x(\ve)$. Now we have
$$\begin{array}{l}
\Prob\{|A(n(\ve))-x(\ve)|<\dt(\ve)\}
\\
=\play{\int_{|\sum_{k=1}^{n(\ve)} (x_k -
x(\ve))|<n(\ve)\dt(\ve)}\Prob(dx_1)...\Prob(dx_{n(\ve)})}
\\
\geq \play{\exp(-n(\ve)\dt(\ve)|\eta|)\exp(-n(\ve)
x(\ve)\eta)\int_{|\sum_{k=1}^{n(\ve)} (x_k -
x(\ve))|<n(\ve)\dt(\ve)}\exp(\eta \sum_{k=1}^{n(\ve)}
x_k)\Prob(dx_1)...\Prob(dx_{n(\ve)})}
\\
=\exp(-n(\ve)\dt(\ve)|\eta|)\exp(-n(\ve)(
x(\ve)\eta-\cH_0(\eta)))\widehat{\Prob}^\ve\{|A(n(\ve))-x(\ve)|<\dt(\ve)\}
\\
\geq
\exp(-n(\ve)\dt(\ve)b)\exp(-n(\ve)\cL_0(x(\ve)))\widehat{\Prob}^\ve\{|A(n(\ve))-x(\ve)|<\dt(\ve)\}
\ .
\end{array}$$

As we have, in this case
$\widehat{\Prob}^\ve\{|A(n(\ve))-x(\ve)|\geq\dt(\ve)\}\leq
\dfrac{\widehat{\E}^\ve|\zt-x(\ve)|^2}{n(\ve)\underline{\dt}^2}\ra
0$ as $\ve \da 0$, we have
$\widehat{\Prob}^\ve\{|A(n(\ve))-x(\ve)|<\dt(\ve)\}\ra 1$ as $\ve
\da 0$. We choose $\ve_0=\ve_0(\nu,\underline{\dt})$ small enough
such that for $0<\ve<\ve_0$ we have
$\widehat{\Prob}^\ve\{|A(n(\ve))-x(\ve)|<\dt(\ve)\}\geq
\exp(-n(\ve)\nu/2)$. We then choose $\overline{\dt}$ small enough
such that $\overline{\dt}b\leq \nu/2$. This then gives the lower
bound. $\square$

\

The next lemma gives some simple but important properties of the
functions $\cH(h,\bt)$ and $\cL(h,\al)$, which will be used later.

Let us denote $\bt[h,\al]=\text{argmax}_{\bt}(\al\bt-\cH(h,\bt))$.
Let $h\in [H_0, \Hbar]$.

\

\textbf{Lemma 3.2.} \textit{We have}

\textit{1. $\cH(h,0)=0$ and $\dfrac{\pt}{\pt
\bt}\cH(h,\bt)|_{\bt=0}<0$.}

\textit{2. For any $b>0$ the functions $\cH$ and $\dfrac{\pt}{\pt
\bt}\cH$ are uniformly continuous in $(h,\bt)$, $|\bt|<b$. The
function $\cH(h,\bt)$ is $C^{\infty}$ in the variables $h$ and
$\bt$.}

\textit{3. The function $\cH(h,\bt)$ is strictly convex in $\bt$.}

\textit{4. We have $\left|\dfrac{\pt}{\pt \bt}\cH(h,\bt)\right|\leq
U\sqrt{2h}$ for some constant $U>0$. When $|\al|>U\sqrt{2h}$ we have
$\cL(h,\al)=+\infty$.}

\textit{5. The set $A(h)=\{\al: \cL(h,\al)<\infty\}$ has nonempty
interior.}

\textit{6. Let $\widehat{\al}$ be such that $\cL(h,
\widehat{\al})=0$, then $\widehat{\al}$ is in the interior of the
set $A(h)$.}

\textit{7. Let $|\bt[h,\al]|\leq b<\infty$. Then for any small
$\kp>0$ and any $|\al'-\al|<\kp$, $|h'-h|<\kp$ we have
$|\bt[h,\al]-\bt[h',\al']|<C(b,\kp)$ and
$|\cL(h,\al)-\cL(h',\al')|<C(b,\kp)$ for a constant $C(b,\kp)\da 0$
as $\kp \da 0$.}

\

\textbf{Proof.} For notational convenience let $\zt=-(\xi+\eta)$.

1. Let $\cH_0(\beta)=\ln \mathbb{E}\exp(\beta\zt)$. We have
$\mathcal{H}(h,\beta)=\play{\frac{\sqrt{2h}}{2D}}\cH_0(\beta)$. It
is obvious that $\cH(h,0)=0$. Also we have $\dfrac{\pt}{\pt
\bt}\cH(h,\bt)|_{\bt=0}=\dfrac{\sqrt{2h}}{2D}\E \zt<0$.

2. We have $$\cH(h,\bt)=\dfrac{\sqrt{2h}}{2D}\ln \E \exp(\bt\zt)$$
and
$$\play{\dfrac{\pt}{\pt \bt}
\cH(h,\bt)=\dfrac{\sqrt{2h}\E\zt\exp(\bt\zt)}{2D \E\exp(\bt\zt)}}
$$ so that they are uniformly continuous in $(h,\bt)$ for $|\bt|<b$.
One can take higher derivatives also so that the function
$\cH(h,\bt)$ is $C^\infty$ in both variables $h$ and $\bt$.

3. We can calculate $$\begin{array}{l} \dfrac{\pt ^2}{\pt
\bt^2}\cH(h,\bt)=\dfrac{\sqrt{2h}(\E\zt^2\exp(\bt \zt)\E\exp(\bt
\zt)-(\E \zt\exp(\bt \zt))^2)}{2D(\E\exp(\bt\zt))^2}>0 \
.\end{array}$$ since the Cauchy-Schwarz inequality is now a strict
one. This means that the function $\cH(h,\bt)$ is strictly convex in
$\bt$.

4. From 2 we can have $\left|\dfrac{\pt}{\pt
\bt}\cH(h,\bt)\right|\leq U \sqrt{2h}$. This gives, when
$|\al|>U\sqrt{2h}$, that $\cL(h,\al)=+\infty$.

5. Since we assumed that $\Prob\{|\xi|\leq M\}=\Prob\{|\eta|\leq
M\}=1$ and both have density, we can assume that there exist $c<0,
C>0$ such that $\Prob(c<\xi+\eta<C)=1$ and there exist $\kp>0$ and
$\mu>0$ such that $\Prob(\xi+\eta>C-\mu)\geq \kp$ and
$\Prob(\xi+\eta<c+\mu)\geq \kp$, also $C-\mu>c+\mu$. From here we
get, that for $\bt>0$, $\cH(h,\bt)\geq
-\dfrac{\sqrt{2h}}{2D}(c+\mu)\kp \bt$ and for $\bt<0$,
$\cH(h,\bt)\geq -\dfrac{\sqrt{2h}}{2D}(C-\mu)\kp \bt$. This fact
helps us to conclude that $\{\al\in \R,
-\dfrac{\sqrt{2h}}{2D}(C-\mu)\kp<\al<
-\dfrac{\sqrt{2h}}{2D}(c+\mu)\kp\}\subset A^\circ(h)$.

6. We are proving that the number $\widehat{\al}$ which makes
$\cL(h, \widehat{\al})=0$ is in the interior of the set $A(h)$.
Since $\widehat{\al}=\dfrac{\pt}{\pt
\bt}\cH(h,\bt[h,\widehat{\al}])$ and
$\widehat{\al}\bt[h,\widehat{\al}]=\cH(h,\bt[h,\widehat{\al}])$. By
strict convexity of $\cH$ in $\bt$ this means that
$\widehat{\al}=\dfrac{\pt}{\pt \bt}\cH(h,\bt)|_{\bt=0}$. The
statement reduces to proving that
$-\dfrac{\sqrt{2h}}{2D}(C-\mu)\kp<\dfrac{\pt}{\pt\bt}\cH(h,\bt)|_{\bt=0}<
-\dfrac{\sqrt{2h}}{2D}(c+\mu)\kp$, which is straightforward.

7. Suppose $|\bt[h,\al]|\leq b<\infty$. Then by strict convexity of
$\cH(h,\bt)$ in $\bt$ we see that $\bt[h,\al]$ is the unique
solution of the equation $\cH_0'(\bt)=\dfrac{2D\al}{\sqrt{2h}}$.
This also gives $|\al|\leq K(b)$ for some constant $K(b)>0$. For any
$|\al'-\al|<\kp$ and $|h'-h|<\kp$ we have
$\play{\left|\dfrac{2D\al}{\sqrt{2 h}}-\dfrac{2D
\al'}{\sqrt{2h'}}\right|< V \kp}$ for some $V>0$. Therefore from the
smoothness of the function $\cH_0(\bt)$ and the strict monotonicity
of $\cH_0'(\bt)$ we conclude that the unique solution $\bt[h',\al']$
of the equation $\cH'_0(\bt)=\dfrac{2D \al'}{\sqrt{2h'}}$ is close
to $\bt[h,\al]$: $|\bt[h,\al]-\bt[h',\al']|<C(b,\kp)$. This gives
also the fact that

$$\begin{array}{l}
|\cL(h,\al)-\cL(h',\al')|
\\
\leq|\al\bt[h,\al]-\al'\bt[h',\al']|+|\cH(h,\bt[h,\al])-\cH(h',\bt[h',\al'])|
\\
\leq|\al||\bt[h,\al]-\bt[h',\al']|+|\al-\al'||\bt[h',\al']| +
\\
\ \ \ \ \ \ \ \ \ \ +|\cH(h,\bt[h,\al])-\cH(h',\bt[h,\al])| +
|\cH(h',\bt[h,\al])-\cH(h',\bt[h',\al'])|
\\
< C(b,\kp)
\end{array}$$
for some $C(b,\kp)>0$, and we have $C(b,\kp)\da 0$ as $\kp \da 0$.
$\square$

\

The next lemma is an analogue of Lemma 5 in [7].

\

\textbf{Lemma 3.3.} \textit{For any $\nu>0$ there exist some
$\Dt(\nu)>0$, $\dt_0(\nu)>0$ such that for any fixed
$0<\dt_0<\dt_0(\nu)$ and fixed $0<\Dt<\Dt(\nu)$, there exist
$\dt_1(\nu,\Dt)>0$ such that for any $0<\dt_1<\dt_1(\nu,\Dt)$ on the
set $|\hH_{t_0}^\ve-h_0|<\dt_0$, uniformly with respect to $t_0 ,
h_0 ,
 H_{t_0}^\ve  , q_0$, under the condition
$\left|\bt[h_0, \dfrac{h_1-h_0}{\Dt}]\right|\leq b<\infty$, as $\ve
\da 0$ , we have}
$$\begin{array}{l}
\exp(-\ve^{-1}(\Dt( \cL(h_0,
\dfrac{h_1-h_0}{\Dt}))-C(b)\nu\Dt-\widetilde{C}(b,\dt_0)))
\\
\ \ \ \ \ \geq \Prob\{|\hH_{t_0+\Dt}^\ve-h_1|<\dt_1|\cF_{t_0}\}\geq
\\
\ \ \ \ \ \ \ \ \ \ \exp(-\ve^{-1}(\Dt( \cL(h_0,
\dfrac{h_1-h_0}{\Dt}))+C(b)\nu\Dt+\widetilde{C}(b,\dt_0)))
\end{array}$$
\textit{where $C(b)>0$ is a constant and $\widetilde{C}(b,\dt_0)\da
0$ as $\dt_0 \da 0$.}

\

\textbf{Proof.} Let $\zt=-(\xi+\eta)$. Let $\cH_0(\bt)=\ln \E \exp
(\bt \zt)$. Let $\cL_0(\al)=\sup\li_{\bt}(\al \bt - \cH_0(\bt))$.
Let $\Dt>0$. Let $N^\ve(t_0,\Dt)$ be the number of crossings that
the process $Q_t^\ve$ make with the set $\{Q=kD, k\in \N\}$ during
time $[t_0, t_0+\Dt]$. Let $n^\ve(t_0,\Dt)=N^\ve(t_0, \Dt)/2$ if
$N^\ve(t_0,\Dt)$ is even and $n^\ve(t_0,\Dt)=(N^\ve(t_0,\Dt)-1)/2$
if $N^\ve(t_0,\Dt)$ is odd. Since
$$\dfrac{1}{\ve}\sqrt{2H_0}\leq\dot{Q}_t^\ve=\dfrac{1}{\ve}\sqrt{2H_t^\ve}\leq
\dfrac{1}{\ve}\sqrt{2\Hbar}$$ we have $$\play{\dfrac{1}{\ve}\Dt
\sqrt{2\Hbar}\geq \int_{t_0}^{t_0+\Dt}\dot{Q}_t^\ve dt\geq
\dfrac{1}{\ve}\Dt \sqrt{2H_0}} \ .$$ This together with the fact
that $$(N^\ve(t_0,\Dt)-1)D\leq
\play{\int_{t_0}^{t_0+\Dt}}\dot{Q}_t^\ve dt\leq N^\ve(t_0,\Dt)D \
.$$ gives $$C_1 \dfrac{\Dt}{\ve}\leq N^\ve(t_0,\Dt)\leq
C_2\dfrac{\Dt}{\ve}$$ for some $C_1>0 , C_2>0$. Since we assumed
that $\Prob\{|\xi|\leq M\}=\Prob\{|\eta|\leq M\}=1$ we see that
$|H_t^\ve-H_{t_0}^\ve|\leq C_3\Dt$ for $t\in [t_0, t_0+\Dt]$ and
some $C_3>0$. As we have $|H_{t_0}^\ve-h_0|<\dt_0$ we have
$|H_t^\ve-h_0|<\dt_0+C_3\Dt$ for $t\in [t_0,t_0+\Dt]$. This gives
$$|\sqrt{2H_t^\ve}-\sqrt{2h_0}|\leq C_4(\dt_0+\Dt)$$ for some $C_4>0$
and $t\in [t_0, t_0+\Dt]$. Therefore
$$\left|N^\ve(t_0,\Dt)D-\dfrac{1}{\ve}\sqrt{2h_0}\Dt\right|\leq
\left|\int_{t_0}^{t_0+\Dt}(\dot{Q}_t^\ve-\dfrac{1}{\ve}\sqrt{2h_0})
dt\right|+D\leq \dfrac{C_4 (\dt_0+\Dt)}{\ve}\Dt+D \ .$$ This gives
$$\left|\ve n^\ve(t_0,\Dt)-\dfrac{\sqrt{2h}}{2D} \Dt\right|\leq
C_5(\ve+(\dt_0+\Dt)\Dt)$$ for some $C_5>0$. This implies that for
$\ve+(\dt_0+\Dt)\Dt<<\Dt$ we have $n^\ve(t_0,\Dt)\ra \infty$ as $\ve
\da 0$. Also, in this case $C_6\Dt\leq\ve n^\ve(t_0,\Dt)\leq C_7
\Dt$ for some $C_6, C_7>0$.

Let $\zt_k=-(\xi_k+\eta_k)$. Now we have, for $\ve>0$ small enough,

$$\begin{array}{l}
\Prob\{|H_{t_0+\Dt}^\ve-h_1|<\dt_1|\cF_{t_0}\}
\\
\play{=\Prob\{|H_{t_0+\Dt}^\ve-H_{t_0}^\ve+H_{t_0}^\ve-h_1|<\dt_1|\cF_{t_0}\}}
\\
\play{\geq \Prob\left\{\left|\ve
(\zt_1+...+\zt_{n^\ve(t_0,\Dt)})+H_{t_0}^\ve-h_1\right|<\dt_1/2
|\cF_{t_0}\right\}} \ .
\end{array}$$

Fix some $\lb>0$ such that $C_5\lb<\dfrac{\sqrt{2h}}{8D}$. We then
choose $\dt_0(\nu)$ and $\Dt(\nu)$ such that
$\dt_0(\nu)+\Dt(\nu)<\lb$ and we fix some $0<\dt_0<\dt_0(\nu)$ and
$0<\Dt<\Dt(\nu)$. We then choose $\ve$ small enough such that
$C_5\ve/\Dt<\dfrac{\sqrt{2h}}{8D}$. We see that for a chosen
$\dt_1(\nu,\Dt)>0$ such that $\dt_1(\nu,\Dt)/\Dt$ is small, for any
$0<\dt_1<\dt_1(\nu,\Dt)$, we can make $\dfrac{\dt_1/2}{\ve
n^\ve(t_0,\Dt)}$ to be smaller than the $\overline{\dt}$ in Lemma
3.1. Also, we notice that $\dfrac{\dt_1/2}{\ve n^\ve(t_0,\Dt)}$ is
bounded away from $0$ as $\ve \da 0$, for fixed $\dt_1$ and $\Dt$.
On the other hand, since we have
$$
\play{\text{argmax}_{\bt}\left(\dfrac{h_1-H_{t_0}^\ve}{\ve
n^\ve(t_0,\Dt)}\bt -\cH_0(\bt)\right)}
\play{=\text{argmax}_{\bt}\left((h_1-H_{t_0}^\ve)\bt -\ve
n^\ve(t_0,\Dt) \cH_0(\bt)\right)}
$$
which, by Lemma 3.2.7, is close to
$$\bt[h_0,\dfrac{h_1-h_0}{\Dt}]=\text{argmax}_\bt((h_1-h_0)\bt-\dfrac{\sqrt{2h_0}}{2D}\Dt\cH_0(\bt)) \ ,$$
say, within a distance of $\kp(\dt_0,(\dt_0+\Dt)\Dt)$, as $\ve$ is
small. And this $\kp(\dt_0, (\dt_0+\Dt)\Dt)\ra 0$ as $(\dt_0,\Dt)\ra
(0,0)$. We shall choose our $\dt_0(\nu)$ and $\Dt(\nu)$ to be small
such that $|\kp(\dt_0(\nu),(\dt_0(\nu)+\Dt(\nu))\Dt(\nu))|<1$. By
our assumption $\left|\bt[h_0,\dfrac{h_1-h_0}{\Dt}]\right|\leq b
<\infty$. Now we see that Lemma 3.1 applies. We then get, on the set
$\{|H_{t_0}^\ve-h_0|<\dt_0\}$, as $\ve$ is small, we have

$$\begin{array}{l}
\Prob\{|H_{t_0+\Dt}^\ve-h_1|<\dt_1|\cF_{t_0}\}
\\
\play{\geq\Prob\left\{\left|\dfrac{\zt_1+...+\zt_{n^\ve(t_0,\Dt)}}{n^\ve(t_0,\Dt)}-\dfrac{h_1-H_{t_0}^\ve}{\ve
n^\ve(t_0,\Dt)}\right|<\dfrac{\dt_1/2}{\ve n^\ve(t_0,\Dt)}
|\cF_{t_0}\right\}}
\\
\geq \exp(- n^\ve(t_0,\Dt)(\sup\li_{\bt\in
\R}\left(\dfrac{h_1-H_{t_0}^\ve}{\ve
n^\ve(t_0,\Dt)}\bt-\cH_0(\bt)\right)+\nu)) \ .
\end{array}$$

Now we use Lemma 3.2.7 to get the bound

$$\begin{array}{l}
\Prob\{|H_{t_0+\Dt}^\ve-h_1|<\dt_1|\cF_{t_0}\}
\\
\geq \exp(- \ve^{-1} \left(\ve n^\ve(t_0,\Dt)\sup\li_{\bt\in
\R}\left(\dfrac{h_1-h_0}{\ve
n^\ve(t_0,\Dt)}\bt-\cH_0(\bt)\right)+C_8\nu\Dt+C_9(b,\dt_0)\right))
\\
= \exp(-\ve^{-1}(\Dt\sup\li_{\bt\in
\R}\left(\dfrac{h_1-h_0}{\Dt}\bt-\dfrac{\ve
n^\ve(t_0,\Dt)}{\Dt}\cH_0(\bt)\right)+C_8\nu\Dt+C_9(b,\dt_0)))
\\
\geq\exp(-\ve^{-1}(\Dt\cL(h_0,\dfrac{h_1-h_0}{\Dt})+C_8\nu\Dt+C_9(b,\dt_0)+C_{10}(b,\dt_0,\Dt)\Dt))
\ .
\end{array}$$

Here the auxiliary constant $C_8>0$ and positive functions
$C_{10}(b,\dt_0,\Dt)\ra 0$ as $(\dt_0,\Dt)\ra (0,0)$ and
$C_9(b,\dt_0)\ra 0$ as $\dt_0 \ra 0$. We choose $\dt_0(\nu)$ and
$\Dt(\nu)$ small enough such that
$C_8\nu+C_{10}(b,\dt_0(\nu),\Dt(\nu))\leq C_{11}(b)\nu$ for
$C_{11}>0$. This gives, for some $C(b)>0$ and
$\widetilde{C}(b,\dt_0)\da 0$ as $\dt_0 \da 0$, the bound

$$\Prob\{|H_{t_0+\Dt}^\ve-h_1|<\dt_1|\cF_{t_0}\} \geq
\exp(-\ve^{-1}(\Dt( \cL(h_0,
\dfrac{h_1-h_0}{\Dt}))+C(b)\nu\Dt+\widetilde{C}(b,\dt_0))) \ .$$

As we have (2.4), this also gives, as $\ve$ is small, that on the
set $\{|\hH_{t_0}^\ve-h_0|<\dt_0\}$ we have

$$\Prob\{|\hH_{t_0+\Dt}^\ve-h_1|<\dt_1|\cF_{t_0}\} \geq
\exp(-\ve^{-1}(\Dt( \cL(h_0,
\dfrac{h_1-h_0}{\Dt}))+C(b)\nu\Dt+\widetilde{C}(b,\dt_0))) \ .$$

Similarly one can estimate

$$\Prob\{|\hH_{t_0+\Dt}^\ve-h_1|<\dt_1|\cF_{t_0}\} \leq
\exp(-\ve^{-1}(\Dt( \cL(h_0,
\dfrac{h_1-h_0}{\Dt}))-C(b)\nu\Dt-\widetilde{C}(b,\dt_0))) \ .$$
 $\square$

\

\textbf{Remark.} In the proof of Theorem 3.1 we will iteratively use
this Lemma and we emphasize that the choice of $\dt_1$ does not
depend on the choice of $\dt_0$ (of course, provided that
$\dt_0<\dt_0(\nu)$). Also, the choice of small $\ve$ may depend on
$\nu$, $\Dt$, $\dt_1$, $\dt_1/\Dt$ (coming from the dependence of
$\ve$ on $\underline{\dt}$ in Lemma 3.1).

\

\textbf{Proof of Theorem 3.1.}

1. \textit{Set-up}. Let $\mathcal{L}$ be the Legendre transform of
$\mathcal{H}$:

$$\mathcal{L}(h,\alpha)=\sup\li_{\beta}(\alpha\beta-\mathcal{H}(h,\beta)).$$

And the action functional is defined as

$$S_{0T}(\varphi)=\left\{
\begin{array}{l}
\play{\int_0^T \mathcal{L}(\varphi_s, \dot{\varphi}_s)ds}  \ ,
\text{ for} \ \varphi\in C_{0T}([H_0,\Hbar]) \ \text{
absolutely \ continuous} \ , \\
+\infty \ ,  \ \text{otherwise} \ .
\end{array}\right.$$

Part (0) of the large deviation principle can be shown as Lemma
7.4.2 of [5].

2. \textit{First part of the proof}. The lower bound (I).

Let $S(\vphi)<\infty$. We show that, given any $\nu>0$, any $\dt>0$,
we have for $\ve$ small enough, that

$$\ve\ln \Prob\{\rho_{0T}(\widehat{H}^\ve, \vphi)<\dt\}\geq -S(\vphi)-\nu \, .$$

Assume that for any $s$, $\cL(\vphi_s, \dot{\vphi}_s)<\infty$ for
any $s$. The reason is the same as in [7], Section 4, Step 1. By
Lemma 3.2.4 we can assume that $\sup\li_{0\leq s\leq
T}|\dot{\vphi}_s|\leq U\sqrt{2\Hbar}$ ($U$ is the constant coming
from Lemma 3.2.4).

3. By Lemma 3.2.5 for each $s\in [0,T]$ the set $\{\al: \cL(\vphi_s,
\al)<\infty\}$ has non-empty interior $\cL^\circ[\vphi_s]$.  Since
$\cL(\vphi_s, \dot{\vphi}_s)<\infty$ we have, as in Section 4, Step
5 of [7] that $\cL(\vphi_s, \dot{\vphi}_s)=\liminf\li_{\al \ra
\dot{\varphi}_s, \al\in \cL^\circ[\varphi_s]}\cL(\varphi_s, \al)$.
For each such $\al$ there exists a (unique in our case since
$\cH(h,\bt)$ is strictly convex in $\bt$) finite adjoint
$\bt[\vphi_s, \al]$. We have $\cH(\vphi_s,\bt)=\al\bt[\vphi_s,
\al]-\cL(\vphi_s,\al)$ and $\cL(\vphi_s,\al)=\al\bt[\vphi_s,
\al]-\cH(\vphi_s,\bt)$. We then choose
$\dot{\widetilde{\vphi}}_s=\al\in \cL^\circ[\vphi_s]$ so that the
value $\cL(\vphi_s, \dot{\widetilde{\vphi}}_s)$ is close to
$\cL(\vphi_s, \dot{\vphi}_s)$.

Put $\widetilde{\vphi}_t=\hH_0^\ve+\play{\int_0^t
\dot{\widetilde{\vphi}}_sds}$. We can choose this new curve to be as
close to $\vphi_t$ as we like, therefore we can make

$$\left|\int_0^T \cL(\vphi_s, \dot{\widetilde{\vphi}}_s)ds-S_{0T}(\vphi)\right|\leq \nu/3 \ .$$

4. For any $s$ we choose a measurable $\widehat{\al}_s$ such that
$\cL(\vphi_s, \widehat{\al}_s)=0$. This is the same as in Section 4,
Step 6 of [7]. We see that $|\widehat{\al}_s|\leq U\sqrt{2\Hbar}$
for the constant $U$ in Lemma 3.2.4. Also, $\widehat{\al}_s\in
\cL^\circ[\vphi_s]$ by Lemma 3.2.6. For this $\widehat{\al}_s$ the
corresponding adjoint $\bt[\vphi_s,
\widehat{\al}_s]=\text{argmax}_\bt(\bt \widehat{\al}_s-\cH(\vphi_s,
\bt))$ exists, is unique and finite.

5. As is the same in Section 4, Step 7 of [7], we take for given $b$
that

$$\vphi_t^b=\hH_0^\ve+
\int_0^t \left(\dot{\widetilde{\vphi}}_s\1(|\bt[\vphi_s,
\dot{\widetilde{\vphi}}_s]|\leq b)+\widehat{\al}_s\1(|\bt(\vphi_s,
\dot{\widetilde{\vphi}}_s)|>b)\right)ds \ .$$

Since we choose $|\widehat{\al}_s|\leq U\sqrt{2\Hbar}$ in Step 4 of
our proof we can find a $b$ such that the curve $\vphi^b$ is still
close to $\vphi$ in $\rho_{0T}$ norm, and the values $\play{\int_0^T
\cL(\vphi_s, \dot{\vphi}^b_s)}ds$ and $\play{\int_0^T \cL(\vphi_s,
\dot{\vphi}_s)}ds$ are close to each other, say
$$\left|\int_0^T\cL(\vphi_s, \dot{\vphi}^b_s)ds-S_{0T}(\vphi)\right|\leq
2\nu/3 \, .$$
 At the same time,
we make $|\bt[\varphi_s, \dot{\vphi}^b_s]|\leq b$. And for $b$ large
enough we make $\{\rho_{0T}(\widehat{H}^\ve, \vphi)<\dt\}\supset
\{\rho_{0T}(\hH^\ve, \vphi^b)<\dt/2\}$.

6. Similarly as in Section 4, Step 9 of [7], we change our functions
$\vphi$ and $\vphi^b$ into a step function $\psi$ and a piecewise
linear function $\chi$ on $[0,T]$ so that first
$\psi_s=\vphi_{[s/\Dt]\Dt}$ and
$\dot{\chi}_s=\dot{\vphi}^b_{[s/\Dt]\Dt}$ and the steplength $\Dt$
of $\psi$, $\chi$ satisfies $\Dt<\Dt(\nu)$ (the value from Lemma
3.3). Secondly,
$$\left|\int_0^T\cL(\psi_s, \dot{\chi}_s)-S_{0T}(\vphi)\right|\leq \nu \, .$$
Thirdly,
$$\{\rho_{0T}(\hH^\ve,\vphi)<\dt\}\supset \{\rho_{0T}(\hH^\ve,\chi)<\dt'\} \, $$
if $\dt'$ is small w.r.t $\dt$. At last, our choice of $\chi$ can
make all the Fenchel-Legendre adjoint to the $\dot{\chi}_s$ variable
$\bt[\psi_{(m-1)\Dt},\dot{\chi}_{(m-1)\Dt+}]=\text{argmax}_{\bt}(\bt
\dot{\chi}_{(m-1)\Dt+}-\cH(\psi_{(m-1)\Dt},\beta))$ uniformly
bounded, $|\bt[\psi_{(m-1)\Dt},\dot{\chi}_{(m-1)\Dt+}]|\leq b$.

To achieve these goals, we need to choose $\dt'<<\dt$ and $\Dt$
small enough such that $\rho_{0T}(\vphi, \chi)<\dt'$ is small and
$\rho_{0T}(\chi,\psi)<\dt'$ ($\psi$ is not in the space
$C_{0T}([H_0,\Hbar])$ but we can still use the distance $\rho_{0T}$)
is small. This $\Dt$ is chosen based on given small $\nu$, $\dt$,
$\dt'$ and the value $b$ found in Step 5.

7. The next step is the same as in Section 4, Step 10 of [7]. Let
$N\Dt=T$. Let $\vphi^\Dt=(\vphi_\Dt, ... , \vphi_{N\Dt})$ and
$(\hH^\ve)^\Dt=(\hH^\ve_{\Dt}, ... , \hH^\ve_{N\Dt})$. First of all
we have
$$\{\rho_{0T}(\hH^\ve, \chi)<\dt'\}\supset \{\rho_{0T}^{\text{discrete}}((\hH^\ve)^\Dt, \chi^\Dt)<\dt''\}$$
if ever $\dt''$ and $\Dt$ are small. The $\dt''$ and $\Dt$ are
chosen based on given small $\dt'$. Here
$\rho_{0T}^\text{discrete}(\psi^\Dt,
\chi^\Dt)=\sup\li_{m}|\psi_{m\Dt}-\vphi_{m\Dt}|$ and the inclusion
comes from the fact that as $\Dt$ is small we have (regarding
$\chi^\Dt$ and $(\hH^\ve)^\Dt$ as step functions also)

$$\rho_{0T}(\hH^\ve, \chi)\leq \rho_{0T}^{\text{discrete}}((\hH^\ve)^\Dt, \chi^\Dt)+
\rho_{0T}((\hH^\ve)^\Dt, \hH^\ve)+\rho_{0T}(\chi^\Dt, \chi)\leq
\dt''+C \Dt<\dt' $$ if ever $\dt''<\dt''(\dt')$ and $\Dt\leq
\Dt(\dt')$.

Then we estimate

$$\Prob\{\rho_{0T}^\text{discrete}(\widehat{H}^\Dt,
\chi^\Dt)<\dt''\}
\geq\E\prod_{m=1}^{N}\1(|\hH_{m\Dt}-\chi_{m\Dt}|<\dt_m''')$$ with
$\dt_1'''<\dt_2'''<...<\dt_{N}'''<\dt''\wedge \dt_0(\nu)\wedge
\dt_1(\nu,\Dt)$ to be chosen later ($\dt_0(\nu)$ and
$\dt_1(\nu,\Dt)$ are from Lemma 3.3).

8. Let us estimate the conditional expectation
$\E(\1(|\hH^\ve_{m\Dt}-\chi_{m\Dt}|<\dt_m''')|\cF_{(m-1)\Dt})$ on
the set $\{|\hH^\ve_{(m-1)\Dt}-\chi_{(m-1)\Dt}|<\dt_{m-1}'''\}$.

We apply Lemma 3.3 and Lemma 3.2.7 on the set
$\{|\hH^\ve_{(m-1)\Dt}-\chi_{(m-1)\Dt}|<\dt_{m-1}'''\}$ to get the
estimate

$$\begin{array}{l}
\E(\1(|\hH^\ve_{m\Dt}-\chi_{m\Dt}|<\dt_m''')|\cF_{(m-1)\Dt})
\\
\geq \exp(-\ve^{-1}(\Dt
\cL(\chi_{(m-1)\Dt},\dot{\chi}_{(m-1)\Dt+})+C(b+\kp(\dt'),\dt_{m-1}''')\nu\Dt+\widetilde{C}(b+\kp(\dt'),\dt_{m-1}''')))
\\
\geq \exp(-\ve^{-1}(\Dt
\cL(\psi_{(m-1)\Dt},\dot{\chi}_{(m-1)\Dt+})+A(b,\dt')\Dt+KC(b,\dt_{m-1}''')\nu\Dt+K\widetilde{C}(b,\dt_{m-1}''')))
\\
\geq \exp(-\ve^{-1}(\Dt
\cL(\psi_{(m-1)\Dt},\dot{\chi}_{(m-1)\Dt+})+A(b)\nu\Dt+K\widetilde{C}(b,\dt_{m-1}'''))
 \ .
\end{array}$$

Here $\kp(\dt')\ra 0$ as $\dt'\da 0$. We have used the fact that
$\rho_{0T}(\chi,\psi)<\dt'$ and Lemma 3.2.7, as well as the fact
that $|\bt[\psi_{(m-1)\Dt},\dot{\chi}_{(m-1)\Dt+}]|\leq b$. We are
choosing $\dt'$ small such that
$C(b+\kp(\dt'),\dt_{m-1}''')<KC(b,\dt_{m-1}''')$ and
$\widetilde{C}(b+\kp(\dt'),\dt_{m-1}''')<K\widetilde{C}(b,\dt_{m-1}''')$
for some $K>0$. The constant $A(b,\dt')\ra 0$ as $\dt' \da 0$. We
are choosing $\dt'$ small enough such that
$A(b,\dt')+KC(b,\dt_{m-1}''')\nu<A(b)\nu$ for some $A(b)>0$. The
constants $C$ and $\widetilde{C}$ are those from the statement of
Lemma 3.3.

9. Now the lower estimate follows from a backward induction: we are
choosing at each step $\widetilde{C}(b,\dt_{m-1}''')$ small enough
compared to
$\widetilde{C}(b,\dt_{m}''')-\widetilde{C}(b,\dt_{m-1}''')$, and we
choose $K\widetilde{C}(b,\dt_{N}''')<\nu$. We have

$$\barr{l}
\Prob\{|\hH_{N\Dt}-\vphi_{N\Dt}|<\dt_N''',...,|\hH_{\Dt}-\vphi_{\Dt}|<\dt_1'''\}
\\
\play{\geq \exp(-\ve^{-1}\Dt\sum\li_{m=1}^N
(\cL(\psi_{(m-1)\Dt},\dot{\chi}_{(m-1)\Dt+})+A(b)\nu)-\ve^{-1}\sum\li_{m=1}^N
K \widetilde{C}(b,\dt_{m-1}''') )}
\\
\play{\geq \exp(-\ve^{-1}(\int_0^T\cL(\psi_s,\dot{\chi}_s)ds+A(b)\nu
T+K \widetilde{C}(b,\dt_{N}''')))}
\\ \geq \exp(-\ve^{-1}(S_{0T}(\vphi)+\nu(A(b)T+2))) \earr$$
as $\ve \da 0$.

(We recall that our choice of parameters has the order
$\nu,\dt\mapsto b\mapsto \dt' \mapsto \Dt, \dt''\mapsto
\dt_N'''\mapsto...\mapsto \dt_1'''\mapsto \ve$.)

10. \textit{Second part of the proof.} The upper bound (II). This
part is similar to that of [7] based on our proof for the lower
bound and we omit it. $\square$

\section{Metastability}

This section is devoted to the description of metastability of
multi-well systems. Due to the stochasticity of the limiting process
at interior vetices of the graph $\Gamma$, the metastability
phenomenon in our case will be metastability of probability
distributions rather than metastability of single states (as in the
classical Freidlin-Wentzell theory , see [5, Ch.6] and compare with
[1]). We will explain below what this is through an example.

We consider generic case when all the width and depth of the wells
are different. Consider two vertices $E_1$ and $E_2$ of our graph
$\Gamma$ (see Fig.1). The vertices $E_1$ and $E_2$ might be exterior
vertices (like $V_1, V_2, V_3, V_4$ in Fig.1) or interior vertices
(like $O_5, O_6, O_7$ in Fig.1). We suppose that the energy levels
corresponding to $E_1$ and $E_2$ are $H_{E_1}$ and $H_{E_2}$, and
$H_{E_2}>H_{E_1}$. Let us first assume that $E_1$ and $E_2$ can be
joined by one edge $I_{N(E_1,E_2)}$. Here $N(E_1, E_2)$ is the
number of the well that has energy level between $H_{E_1}$ and
$H_{E_2}$ (recall that under our convention every well has a highest
and lowest energy level). Recall that the well $N(E_1,E_2)$ has
width $D_{N(E_1,E_2)}$. Let

$$V(E_1,E_2)=\inf\{S_{0T}^{N(E_1,E_2)}(\varphi): H_{E_1}\leq\varphi_t\leq H_{E_2} , 0\leq t \leq T<\infty, \varphi_0=H_{E_1}, \varphi_T=H_{E_2}\}  \ ,
\eqno(4.1)$$ and $$V(E_2,E_1)=0 \, . \eqno(4.2)$$

Here the functional $S_{0T}^{N(E_1,E_2)}$ is defined by
$$S_{0T}^{N(E_1,E_2)}(\varphi)=\play{\int_0^T
\mathcal{L}^{N(E_1,E_2)}(\varphi_t, \dot{\varphi}_t)dt}$$ if
$\varphi$ is absolutely continuous and $+\infty$ otherwise.

The function
$$\mathcal{L}^{N(E_1,E_2)}(h,\alpha)=\sup\li_\beta(\alpha\beta-\mathcal{H}^{N(E_1,E_2)}(h,\beta)) \ , \eqno(4.3)$$
where $\alpha, \beta \in \mathbb{R}$ and $h\geq 0$ is the Legendre
transform of the function
$$\mathcal{H}^{N(E_1,E_2)}(h,\beta)=\play{\frac{\sqrt{2h}}{2D_{N(E_1,E_2)}}\ln
\mathbb{E} \exp(-\beta(\xi^{(N(E_1,E_2))}+\eta^{(N(E_1,E_2))}))} \ .
\eqno(4.4)$$ (If we are at some "small" well, i.e., it contains no
smaller wells we make
$\mathcal{H}^{N(E_1,E_2)}(H_{E_1},\beta)=\play{\frac{\sqrt{2H_{E_1}}}{2D_{N(E_1,E_2)}}\ln
\mathbb{E}
\exp(-\beta(\xi^{(N(E_1,E_2))}+\eta^{(N(E_1,E_2))}))\1(\xi^{(N(E_1,E_2))}<0,
\eta^{(N(E_1,E_2))}<0)} \ . $)

In particular, we see that our function $V(E_1,E_2)$ depends on the
width $D_{N(E_1,E_2)}$ of the $N(E_1,E_2)$-th well, the energy
levels $H_{E_1}$ and $H_{E_2}$ of the $N(E_1,E_2)$-th well and
properties of the random variables $\xi^{(N(E_1,E_2))}$ and
$\eta^{(N(E_1,E_2))}$ which give perturbations at the left and right
walls when the particle is in the $N(E_1,E_2)$-th well.

One can verify that $V(E_1,E_2)$ and $V(E_2,E_1)$ define the
"quasi-potential" for all adjacent vertices $E_1$ and $E_2$ (with
$H_{E_2}>H_{E_1}$) on our graph $\Gamma$. To do this, we shall
notice that by similar arguments as we did in Section 3, the action
functional for the perturbed dynamical system
$\widehat{Y}_t^\ve=(\widehat{H}_t^\ve,K(\widehat{H}_t^\ve,
q_t^\ve))$ on the graph $\Gm$ shall be defined by $S_{0T}(\vphi,
K)=\play{\int_0^T \cL(\vphi_s, K(s), \dot{\vphi}_s)ds}$ where
$\cL(\vphi_s, K(s)
,\dot{\vphi}_s)=\sup\li_{\bt}(\dot{\vphi}_s\bt-\cH(\vphi_s,K(s),
\bt))$. Here $\vphi: [0,T]\ra [H_0, \overline{H}]$ is absolutely
continuous (otherwise the action functional is $+\infty$ and
$\overline{H}>\widehat{H}_t^\ve$ for $0\leq t \leq T$). The function
$K(s):[0,T]\ra \{1,2,...,N\}$ where $N$ is the number of edges of
graph $\Gm$. The function $$\cH(h,K,\bt)=\dfrac{\sqrt{2h}}{2D_K}\ln
\E \exp(-\bt(\xi^{(K)}+\eta^{(K)}))$$ whenever $(h,K)$ does not
correspond to the bottom of a "small" well and it is
$\cH(h_0,K_0,\bt)=\dfrac{\sqrt{2h_0}}{2D_{K_0}}\ln \E
\exp(-\bt(\xi^{(K_0)}+\eta^{(K_0)}))\1(\xi^{(K_0)}<0,
\eta^{(K_0)}<0)$ when $(h_0,K_0)$ corresponds to the bottom of a
"small" well. Since we assume that $\E(\xi+\eta)>0$, we find
(compare with the example given in Section 3) that the minimum in
the definition of the quasi-potential between $E_1$ and $E_2$ is
achieved within the class of functions that satisfy $H_{E_1}\leq
\vphi_t\leq H_{E_2}$ for $0\leq t \leq T$, as defined in (4.1).

Now for any two vertices $F_1$ and $F_2$ on the graph $\Gamma$, let
$$V(F_1,F_2)=\min\li_{(E_1,...,E_m)}\sum\li_{i=1}^{m-1} V(E_i,E_{i+1}) \ . \eqno(4.5)$$

Here $(E_1,...,E_m)$ is a path of $\Gamma$ for which $E_1=F_1,
E_m=F_2$ and each pair $E_i$ , $E_{i+1}$ can be joined by an edge of
the graph $\Gamma$.

One can verify that the function $V(F_1,F_2)$ defines the
"quasi-potential" between $F_1$ and $F_2$ , as was defined in [5,
Ch.6].

In particular, one can easily check that for any interior vertex
$O_l$, there is an exterior vertex $V_k$ such that $V(O_l,V_k)=0$.
Therefore interior vertices are \textit{unstable} (compare with [5,
Ch.6, Lemma 6.4.3]).

Now let us consider the example given in Fig.1. We suppose that,
after using (4.1)-(4.4), we have the following: $V(V_1, O_5)=2,
V(V_2, O_5)=1, V(O_5,O_6)=1, V(O_6, O_7)=1, V(V_3,O_6)=6,
V(V_4,O_7)=5 \, $ and $V(O_7,O_6)=V(O_6,O_5)=V(O_6,
V_3)=V(O_5,V_1)=V(O_5, V_2)=V(O_7, V_4)=0$.

Suppose our process $\widehat{Y}_t^\ve=(\widehat{H}_t^\ve, K_t^\ve)$
starts from a point $(H_0, 7)$ with $H_0$ large enough. Here the
process $\widehat{H}_t^\ve$ is the piecewise linear modification
defined at the beginning of Section 2 and
$\widehat{Y}_t^\ve=Y(\widehat{H}_t^\ve, q_t^\ve)$ is the
identification map introduced in Section 1.

Let $Y_t$ be the (weak) limiting process of $\widehat{Y}_t^\ve$ as
$\ve \da 0$ on the graph $\Gamma$. It is a Markov process on
$\Gamma$ which is a deterministic motion within each edge and only
has stochasticity at the interior vertices (see Theorem 2.1 and
Theorem 2.2). In particular, let the branching probabilities at
vertex $O_7$ be given by $p_6$ (for entering $I_6$) and $p_4=1-p_6$
(for entering $I_4$); the branching probabilities at vertex $O_6$ be
given by $p_5$ (for entering $I_5$) and $p_3=1-p_5$ (for entering
$I_3$); and the branching probabilities at vertex $O_5$ be given by
$p_1$ (for entering $I_1$) and $p_2=1-p_1$ (for entering $I_2$).

After long enough finite time, as $\ve\da 0$, the position of the
process $\widehat{Y}_t^\ve$ will be given by a probability
distribution which is approximately $(p_1p_5p_6, p_2p_5p_6, p_3p_6,
p_4)=(p_1p_5p_6, (1-p_1)p_5p_6, (1-p_5)p_6, 1-p_6)$ among the
exterior vertices $(V_1, V_2, V_3, V_4)$. Let us denote the
distribution $(p_1p_5p_6, p_2p_5p_6, p_3p_6, p_4)$ by $U_0$.

Let us now consider behavior of the process $\widehat{Y}_t^\ve$ at
exponentially long time scale $t=t(\ve)\asymp \exp(C\ve^{-1})$. To
this end we first remind the reader of some classical results in [5,
Ch.6]. Consider a set of $K_i$'s, $i=1,...,l$, which are
equilibriums of a deterministic dynamical system, say $Z_t$. Suppose
the corresponding stochastic dynamical system $Z_t^\ve$, which is a
small random perturbation of $Z_t$, satisfies a large deviation
principle with normalizing factor $\dfrac{1}{\ve}$ and the
quasi-potentials between $K_i$ and $K_j$ are given by $V(K_i,K_j)$.
We decompose the set of $K_i$'s into hierarchy of cycles $\pi^{(0)},
\pi^{(1)},..., \pi^{(s)}$, unified into the last cycle of maximal
rank. For any cycle $\pi^{(k)}$, $0\leq k \leq s$, we define

$$C(\pi^{(k)})=A(\pi^{(k)})
-\min\li_{i\in \pi^{(k)}}\min\li_{g\in
G_{\pi^{(k)}}\{i\}}\sum\li_{(m\ra n)\in g}V(K_m, K_n) \, ,
\eqno(4.6)$$ where $$A(\pi^{(k)})=\min\li_{g\in G(L \backslash
\pi^{(k)})}\sum\li_{(m\ra n)\in g}V(K_m, K_n) \, . \eqno(4.7)$$

Here $L$ is the set of indices for the points $K_1,...,K_l$. The set
$G(L \backslash \pi^{(k)})$ is the collection of all $L \backslash
\pi^{(k)}$-graphs and the set $G_{\pi^{(k)}}\{i\}$ is the collection
of all $i$-graphs restricted to $\pi^{(k)}$ (see [5, Ch.6, Section
6]).

Then for sufficiently small $\rho>0$ we have, $$\lim\li_{\ve \ra
0}\ve\ln \mathbb{E}_x^\ve \tau_{\pi^{(k)}}=C(\pi^{(k)})\, ,
\eqno(4.8)$$ uniformly in $x$ belonging to some $\rho$-neighborhood
of the set $\bigcup_{i\in \pi^{(k)}}K_i$, where $\tau_{\pi^{(k)}}$
is the first exit time for the system $Z_t^\ve$ to exit from
$\pi^{(k)}$ (see Theorem 6.6.2 of [5, Ch.6]).

Also, the asymptotic as $\ve \da 0$ exit position
$Z_{\tau_{\pi^{(k)}}}^\ve$ in $L\backslash \pi^{(k)}$ for the system
to exit from $\pi^{(k)}$ is
 given by one of the $K_i$'s which is the end of the chains of arrows in an $L\backslash \pi^{(k)}$ graph
 that minimizes the sum in
 (4.7) (see Theorem 6.6.1 of [5, Ch.6]).

Now let us turn back to our example. We start from the distribution
$U_0=(p_1p_5p_6, p_2p_5p_6, p_3p_6, p_4)$. The cycles of rank $0$
are just the vertices $V_1, V_2, V_3, V_4$ and we call them
$\pi^{(i-1)}=\{V_i\} \, , i=1,...,4$. We calculate $C(\pi^{(0)})=2$,
$C(\pi^{(1)})=1$, $C(\pi^{(2)})=6$, $C(\pi^{(3)})=5$. Therefore by
using (4.8), we see that at time scale
$t=t(\ve)\asymp\exp(\ve^{-1})$, our system will be jumping out from
$V_2$ first. By determining the $\{V_1,V_3,V_4,O_5,O_6,O_7\}$-graph
minimizing the sum in $A(\pi^{(1)})$, the first vertex that it
approaches will be $V_1$ (to be precise, it will be $O_5$ but $O_5$
is unstable). Taking into account that there is a branching
probability at vertex $O_5$, we see that one such transition will
make the distribution be $(p_1p_5p_6+p_1p_2p_5p_6, p_2^2p_5p_6,
p_3p_6, p_4)$and $n$ times such transitions will make the
distribution be
$(p_1p_5p_6+p_1p_2p_5p_6+p_1p_2^2p_5p_6+...+p_1p_2^np_5p_6,
p_2^{n+1}p_5p_6, p_3p_6, p_4)$. Therefore after many such
transitions, when $n$ is very large, the distribution will be
approximately $U_1=(p_5p_6,0,p_3p_6,p_4)$. The distribution $U_1$
will be the "metastable distribution" over time scale
$t=t(\ve)\asymp\exp(\ve^{-1})$ (compare with [1, Theorems 4.1 and
4.2]).

We increase our time scale. Since $C(\pi^{(0)})=2$,
$C(\pi^{(2)})=6$, $C(\pi^{(3)})=5$, at time scale
$t=t(\ve)\asymp\exp(2\ve^{-1})$, the system begins to jump out from
$V_1$ and transit to $V_2$, which makes the distribution
$U_1=(p_5p_6,0,p_3p_6,p_4)$ be $(p_1p_5p_6,p_2p_5p_6,p_3p_6,p_4)$,
$(p_1^2p_5p_6,p_2p_5p_6+p_2p_1p_5p_6,p_3p_6,p_4)$, ...
,$(p_1^{n+1}p_5p_6,p_2p_5p_6+p_2p_1p_5p_6+...+p_2p_1^np_5p_6,p_3p_6,p_4)$,
and so on. But notice that one such transition happens at time scale
$t=t(\ve)\asymp\exp(2\ve^{-1})$, within which transitions from $V_2$
to $V_1$, as described in the above paragraph, will happen many
times. Therefore over time scale $t=t(\ve)\asymp\exp(2\ve^{-1})$,
$U_1=(p_5p_6,0,p_3p_6,p_4)$ will still be the metastable
distribution of our system.

Over time scale $t=t(\ve)\asymp\exp(2\ve^{-1})$, our system has
already formed a cycle $\{V_1,O_5,V_2\}$, which we call $\pi^{(4)}$.
We calculate $C(\pi^{(4)})=3<C(\pi^{(3)})\wedge C(\pi^{(2)})$. That
means, at time scale $t=t(\ve)\asymp\exp(3\ve^{-1})$, jumping out
from cycle $\pi^{(4)}$ happens first. By determining the
$\{V_3,V_4,O_6,O_7\}$-graph minimizing the sum in $A(\pi^{(4)})$, we
will first jump to $V_3$ (again, it is actually $O_6$ but $O_6$ is
unstable). Taking into account of the branching probabilities, this
will make the distribution be finally $U_2=(0,0,p_6,p_4)$. The
distribution $U_2$ is the metastable distribution over time scale
$t=t(\ve)\asymp\exp(3\ve^{-1})$.

We now consider the cycle $\{V_1, O_5, V_2, O_6, V_3\}$ and we call
it $\pi^{(5)}$. We calculate $C(\pi^{(5)})=7$. Since
$C(\pi^{(2)})=6$, $C(\pi^{(3)})=5$, by the same reasoning above,
over time scale $t=t(\ve)\asymp\exp(5\ve^{-1})$, transition from
$V_4$ to $V_3$ happens first and that leads to a new metastable
distribution $U_3=(0,0,1,0)$.

Over time scale $t=t(\ve)\asymp\exp(6\ve^{-1})$, transition from
$V_3$ to $V_4$ happens. By the same reasoning above, we see that
this leads to the fact that the metastable distribution over time
scale $t=t(\ve)\asymp\exp(6\ve^{-1})$ is still $U_3$. After that
time scale, although new transition might still happen, the
metastable distribution will remain to be $U_3$.

\section*{Acknowledgements}

I thank my advisor M.Freidlin for posing this problem to me and for
pointing out the references [1], [3], [6], [7], as well as many
useful discussions.

\section*{References}

[1] A.Athreya, M.Freidlin, Metastability for random perturbations of
nearly-Hamiltonian systems, \textit{Stochastics and Dynamics}, Vol.
\textbf{8}, No.1 (2008) 1-21.

[2] A.Dembo, O.Zeitouni, \textit{Large Deviations techniques and
applications}, Springer, 1998.

[3] W.Feller, \textit{An Introduction to Probability Theory and Its
Applications}, Vol.2, Second Edition, John Wiley and Sons, 1971.

[4] M.Freidlin, W.Hu, On stochasticity in nearly-elastic systems,
\textit{Stochastics and Dynamics}, online,
DOI:10.1142/S0219493711500201.

[5] M.Freidlin, A.Wentzell, \textit{Random perturbations of
dynamical systems}, Springer, 1998.

[6] A.Yu.Veretennikov, On large deviations in the averaging
principle for SDE's with a "full dependence", \textit{Annals of
Probability}, \textbf{27},1999, No.1, 284-296.

[7] A.Yu.Veretennikov, On large deviations in the averaging
principle for SDE's with a "full dependence", correction, arxiv.
MATH.PR.0502098.

\end{document}